\documentclass[11pt,letterpaper,reqno]{amsart}
\usepackage[letterpaper,textwidth=6.45in,textheight=9in,centering]{geometry}
\usepackage{amssymb,booktabs}
\usepackage{tikz}
\usepackage[hidelinks]{hyperref}
\hypersetup{pdftitle={Even quartets and the maximum density of a simplicial 3-cocycle},pdfauthor={John Mackey}}
\newtheorem{theorem}{Theorem}[section]
\newtheorem{proposition}[theorem]{Proposition}
\newtheorem{lemma}[theorem]{Lemma}
\newtheorem{corollary}[theorem]{Corollary}
\theoremstyle{remark}

\newcommand{\F}{\mathbb F_2}
\newcommand{\E}{\mathbb E}

\title[Even quartets and simplicial cocycles]{Even quartets and the maximum density\protect\\ of a simplicial $3$-cocycle}
\author{John Mackey}
\address{Computer Science, Carnegie Mellon University}
\date{}
\begin{document}
\begin{abstract}
Colour the triples of an $n$-element set red or blue, and call a four-set even if it contains an even number of red triples. We prove that the minimum number of even four-sets is $\frac5{16}\binom n4-\Theta(n^3)$. In particular, the maximum asymptotic density of a simplicial $3$-cocycle over the field of two elements is $11/16$. The classical two-part matrix construction attains the limiting bound. We prove that it also determines the unique limiting distribution of quartet parities on every fixed number of vertices. Writing $m_n$ for the minimum density, we bound the limit inferior and limit superior of $n(5/16-m_n)$ between $0.75$ and $0.752723$. The lower bound on even quartets follows from an exact rational inequality on eight vertices. A perturbation of the construction explains all $72$ cases of equality in this finite inequality, of which only $22$ survive in an extremal limit.
\end{abstract}
\maketitle

\section{Introduction}\label{sec:intro}

Colour each triple of a finite set $V$ red or blue. A \emph{quartet} is a four-element subset of $V$. Call it \emph{even} if it contains zero, two or four red triples, and \emph{odd} otherwise. On any five vertices there are an even number of odd quartets, since each red triple is counted twice. In particular, at least one of the five quartets is even. How few quartets can be even simultaneously?

Write $G=(V,\chi)$ for a colouring, where $\chi\colon\binom V3\to\F$ takes the value $1$ on red triples. Its parity function $q=q_\chi\colon\binom V4\to\F$ and even density $d(G)$ are
\[
 q(Q)=\sum_{T\in\binom Q3}\chi(T),\qquad
 d(G)=\frac{\bigl|\{Q\in\binom V4:q(Q)=0\}\bigr|}{\binom{|V|}4}.
\]
For $n\ge4$ let $TC(n)$ be the minimum number of even quartets over all colourings on $n$ vertices, and put $m_n=TC(n)/\binom n4$.

\begin{theorem}\label{thm:main}
The minimum density of even quartets tends to $5/16$. More precisely, there are absolute constants $c_1,c_2>0$ such that, for all sufficiently large $n$,
\[
 \frac{c_1}{n}\ \le\ \frac5{16}-m_n\ \le\ \frac{c_2}{n}.
\]
Thus $TC(n)=\frac5{16}\binom n4-\Theta(n^3)$.
\end{theorem}

There is also a uniqueness statement. Independently assign each vertex to one of two parts with equal probability, and assign an independent fair sign to each cross pair. Pure quartets are even and quartets split three to one are odd. A quartet split two to two is even precisely when the four cross signs on its rectangle have product $+1$. We call the resulting distribution the \emph{fair two-part law}. Section~\ref{sec:construction} gives a triple colouring realizing it.

\begin{theorem}\label{thm:main-law}
Let $G_n$ be triple colourings with $|V(G_n)|\to\infty$ and $d(G_n)\to 5/16$. For every fixed $s\ge4$, the joint distribution of the $\binom s4$ quartet parities on a uniform ordered sample of $s$ distinct vertices converges to the fair two-part law.
\end{theorem}

\subsection{Cocycles and the Tur\'an problem}
The five-set identity is a complete description of the possible parity patterns. Indeed, let $q\colon\binom V4\to\F$ satisfy
\begin{equation}
        \sum_{Q\in\binom W4}q(Q)=0\qquad(W\in\tbinom V5).
        \label{eq:cocycle}
\end{equation}
Fix $v_0\in V$, give all triples through $v_0$ colour zero, and set $\chi(ijk)=q(v_0ijk)$ on every other triple. Here and below, a string of distinct vertices denotes their underlying set. Equation~\eqref{eq:cocycle} gives $q_\chi=q$ also on the quartets avoiding $v_0$. In simplicial terminology, \eqref{eq:cocycle} says that $q$ is a $3$-cocycle on the full simplex, and the construction writes it as a coboundary. Thus, in Kalai's notation for the maximum support size of such a cocycle,
\[
                   f(n,4)=\binom n4-TC(n).
\]
Theorem~\ref{thm:main} gives $f(n,4)=(11/16)\binom n4+\Theta(n^3)$. Kalai~\cite[Sec.~7.1]{Kalai} describes the matrix construction giving the density $11/16$ and records Peled's flag-algebra upper bound of $0.6916+o(1)$ on the maximum density.

McKay and Radziszowski~\cite{McKayRadziszowski} use $TC(n)$ for the minimum size of a quartet family meeting every five-set in an odd number of quartets. By~\eqref{eq:cocycle}, these are exactly the families of even quartets. The corresponding Tur\'an number $T(n,5,4)$ is the minimum size of a quartet family meeting every five-set at least once. Hence $T(n,5,4)\le TC(n)$. The two-part construction for this covering problem appears in work of de Caen, Kreher and Wiseman~\cite{deCaen} and Giraud~\cite{Giraud}. McKay and Radziszowski record the conjecture that $TC(n)=T(n,5,4)$ at every order. Whether a general covering can improve on the asymptotic density $5/16$ remains open~\cite{Sidorenko}.

The parities are the natural objects for the uniqueness assertion. Switching a pair means reversing the colours of all triples containing it. Every pair switch preserves every quartet parity. Conversely, if two colourings induce the same parities, take their difference $h\colon\binom V3\to\F$. Switch each pair $ij$ avoiding $v_0$ for which $h(v_0ij)=1$. These switches realize $h$ on triples containing $v_0$, and the zero quartet sums of $h$ force agreement on every remaining triple. Thus the parities distinguish colourings exactly up to pair switches.

\subsection{Outline of the proof}
On eight vertices there are $893{,}924$ quartet-parity patterns up to isomorphism. We give an inequality for each of them. When averaged in a large colouring, its right-hand side is a sum of nonnegative conditional quadratic forms in the sense of flag algebras~\cite{Razborov}. The certificate has rational coefficients, and every finite comparison is exact. Section~\ref{sec:perturbation} explains an apparent obstacle to identifying equality cases. Perturbing the two-part construction changes its eight-vertex distribution to first order but its objective only to second order. A sharp pointwise inequality must therefore vanish on patterns additional to the ones that appear in the construction. Section~\ref{sec:equality} proves that the extra patterns disappear from every extremal limit, and recovers a single two-part representation on large samples. Section~\ref{sec:finite} gives an explicit lower bound at every order $n\ge8$. Using the equality law, Section~\ref{sec:second-order} places the first finite correction $n(5/16-m_n)$ between $0.75-o(1)$ and $0.752723+o(1)$. Section~\ref{sec:book} proves a quadratic stability inequality for a six-vertex configuration.

\section{The two-part construction}\label{sec:construction}

Let $V=A\sqcup B$ and let $Z\colon A\times B\to\{-1,+1\}$ be a sign matrix. Represent blue by $+1$ and red by $-1$. Give pure triples sign $+1$, give triples $abb'$ sign $-1$, and give $aa'b$ sign $-Z_{ab}Z_{a'b}$, with letters taken from their indicated parts. The product of the four triple signs on a quartet determines its parity. Direct multiplication gives
\begin{center}
\begin{tabular}{@{}lll@{}}
\toprule
Part sizes on the quartet & Product of triple signs & Even probability in the fair law\\
\midrule
$4+0$ or $0+4$ & $+1$ & $1$\\
$3+1$ or $1+3$ & $-1$ & $0$\\
$2+2$, vertices $i,j,x,y$ & $Z_{ix}Z_{iy}Z_{jx}Z_{jy}$ & $1/2$\\
\bottomrule
\end{tabular}
\end{center}
With independent fair labels and cross signs, the probabilities of these three rows are $1/8$, $1/2$ and $3/8$. The expected even density is therefore $1/8+(3/8)/2=5/16$ at every order. In particular $m_n\le5/16$. The density of each fixed sample pattern converges in probability to its expectation. Indeed, indicators on disjoint vertex sets are independent, so the variance is $O(1/n)$.

\begin{proposition}\label{prop:hadamard}
Suppose $|A|=|B|=m\ge2$. Write $r_{ij}=\sum_{x\in B}Z_{ix}Z_{jx}$. The number of even quartets is
\[
  2\binom m4+\binom m2\frac{m(m-2)}4
                +\frac14\sum_{\{i,j\}\in\binom A2}r_{ij}^2.
\]
When the rows of $Z$ are orthogonal, the even density is
\begin{equation}
   d_{2m}=\frac{(m-2)(5m-6)}{4(2m-1)(2m-3)}
         =\frac5{16}-\frac{24m-33}{16(2m-1)(2m-3)}.
   \label{eq:hadamard}
\end{equation}
\end{proposition}
\begin{proof}
For a fixed pair of rows, the signed sum of the rectangle products over column pairs is $(r_{ij}^2-m)/2$. The number of positive rectangles is one half the number of rectangles plus one half this signed sum. Adding over row pairs and counting the pure quartets proves both formulas.
\end{proof}

For $m=2^j$, index the rows and columns by $\F^j$ and take $Z_{ax}=(-1)^{a\cdot x}$. Distinct rows are orthogonal by summation of a nontrivial character. If $N=2m$ is the least such power of two at least $n$, restriction to a uniform $n$-subset preserves the expected density. Thus $m_n\le d_N$, with $n\le N<2n$. Equation~\eqref{eq:hadamard} gives $5/16-d_N=(3/4+o(1))/N$. This proves the lower bound on the finite deficit in Theorem~\ref{thm:main}. Its upper bound follows from Proposition~\ref{prop:finite-collision}.

\section{The first-order constraints on a sharp bound}\label{sec:perturbation}

Start with the fair two-part construction and reverse each pure $A$-triple independently with probability $t$. Write $d(t)$ for the resulting even probability. Only two kinds of quartet contribute to its change. A quartet in $A$ stays even with probability $(1+(1-2t)^4)/2$. A quartet with three vertices in $A$ and one in $B$ becomes even with probability $t$. Consequently
\begin{align}
 d(t)-\tfrac5{16}
  &=\tfrac1{16}\bigl(\tfrac{1+(1-2t)^4}{2}-1\bigr)+\tfrac14t
    =\tfrac34t^2-t^3+\tfrac12t^4 . \label{eq:perturb}
\end{align}
In particular, the two linear contributions cancel.

This cancellation determines which equalities a finite certificate must allow. Suppose its pointwise remainder is a fixed nonnegative function $R$ of the eight-vertex pattern $H$, and averaging gives $d-5/16\ge\E R(H)$ for every limiting law. Each pattern of positive probability in the fair law must have $R(H)=0$. Now fix a particular eight-vertex construction and a specified triple in $A$. The probability of its prescribed part labels and cross signs is a positive constant. Conditional on these choices, the probability that just the specified triple is flipped is $t(1-t)^{\binom{|A|}3-1}$. A positive remainder on the resulting pattern would therefore force a positive linear term in~\eqref{eq:perturb}. The same conclusion holds for a flip in $B$ by symmetry of the unperturbed quartet law. Up to relabelling there are $22$ original patterns and $50$ further patterns obtainable by one same-side flip. The exact enumeration in Lemma~\ref{lem:certificate} shows that our remainder vanishes on precisely this union.

\section{The eight-vertex inequality}
\label{sec:certificate}

Write $H$ for an unlabelled parity pattern on eight vertices and $e(H)$ for the number of its even quartets. If $H$ is induced by a uniformly sampled eight-set in a colouring $G$, then $\mathbb E[e(H)/70]=d(G)$. Our lower bound comes from the decomposition
\begin{equation}
  \frac{e(H)}{70}-\frac5{16}=\Phi(H)+R(H).
  \label{eq:host-residual}
\end{equation}
The two summands have different positivity properties. The residual $R(H)$ is nonnegative on each eight-vertex pattern. The term $\Phi(H)$ is a positive combination of rooted products, each of which has nonnegative expectation in the limit of growing colourings.

\subsection{Rooted products}

For $r\in\{0,2,4,6\}$ put $k=(8-r)/2$. Choose a uniform ordered $r$-tuple of distinct vertices of $H$ as a root. Split its complement uniformly into two $k$-sets $U,W$. A \emph{flag} is the parity pattern induced by the root and one of the two sets, with root labels retained and vertices in the set unlabelled. Denote the flags by $I,J$. Each term specifies a labelled root pattern $\sigma$ and a symmetric matrix $K$ indexed by flags extending $\sigma$. Its value is defined to be zero at all other root patterns. With $\kappa_0=\kappa_2=\kappa_4=1$ and $\kappa_6=6!$, its contribution before multiplication by its rational weight is
\begin{equation}
  T_{\sigma,K}(H)=
  \kappa_r\,\mathbb E\bigl[\mathbf 1_{\{\mathrm{root}=\sigma\}}K(I,J)\bigr].
  \label{eq:rooted-product}
\end{equation}
The expectation is over the root and split in this one eight-vertex pattern. The symmetry of $K$ makes the order of $U,W$ immaterial.

For an ordinary square the certificate gives an integer function $u\colon\mathcal F_\sigma\to\mathbb Z$ on the set $\mathcal F_\sigma$ of extension flags, and $K(i,j)=u(i)u(j)$. The remaining matrices are specified by a triangle-free graph on a subset of the flags and nonnegative integer amplitudes $a_i$. Put $a_i=0$ off that subset. For adjacency matrix $A$ the entries are $K(i,j)=a_i a_j(1-2A_{ij})$. This is a copositive kernel, meaning $x^{\mathsf T}Kx\ge0$ for every nonnegative vector $x$. Indeed, absorb the amplitudes into the coordinates and normalize the resulting vector $y$ to have sum one. It suffices to prove $\sum_{\{i,j\}\in E(A)}y_i y_j\le1/4$, with each edge counted once. If two nonadjacent vertices both carry mass, the edge sum is affine as their combined mass is divided between them. Moving their combined mass to one endpoint preserves or increases the sum. Iterating leaves a clique, and a clique in a triangle-free graph has at most two vertices. Their product is at most $1/4$. This is the triangle-free case of the Motzkin--Straus theorem~\cite{MotzkinStraus}.

The full sum is $\Phi(H)=\sum_\ell w_\ell T_{\sigma_\ell,K_\ell}(H)$ with $w_\ell>0$. Table~\ref{tab:kernel-counts} gives the 841 terms. The last column gives the divisor after summing the indicated root orders and unordered splits.
\begin{table}[htb]\centering
\begin{tabular}{rrrrr}
\toprule
$r$ & $k$ & Ordinary squares & Triangle-free kernels & Divisor\\
\midrule
0 & 4 & 1   & 0  & 35\\
2 & 3 & 4   & 0  & 560\\
4 & 2 & 116 & 12 & 5040\\
6 & 1 & 707 & 1  & 28\\
\bottomrule
\end{tabular}
\caption{Rooted products in the certificate and their normalization.}\label{tab:kernel-counts}\end{table}
At the empty root the function is $11$ on an even quartet and $-5$ on an odd quartet. At a two-vertex root, a five-set flag is recorded by the quartet bits obtained by omitting either root vertex, followed by the number of odd quartets obtained by omitting an outside vertex. The five-set parity equation leaves eight possible records. There are $28$ root pairs, two orders and ten complementary triple splits. At a four-vertex root, each of the two root parities has 288 extension flags. There are $70$ roots, $24$ orders and three complementary pair splits. At a six-vertex root, the certificate fixes one labelled representative for each of the $16$ root types. Its one-vertex extension has ten free bits. For each compatible six-subset we sum the products over all isomorphisms to the specified representative, then divide by the $28$ possible six-subsets. This sum over root isomorphisms accounts for $\kappa_6$ in~\eqref{eq:rooted-product}.

If $G$ has $n$ vertices, the distribution obtained by first sampling its eight-set and then a root and split is also obtained by sampling the root in $G$ and two disjoint $k$-sets outside it. Instead take two independent uniform $k$-sets outside the root. Conditional on the root, let $x_i$ be the probability of flag $i$ on one extension. For an ordinary term the conditional product expectation is $(\sum_i u(i)x_i)^2$. For a general term it is $x^{\mathsf T}Kx\ge0$. The probability that the two independent sets intersect is $O(1/n)$ and the kernels are fixed and bounded. Conditioning on disjointness therefore changes the expectation by $O(1/n)$, uniformly in $G$. In particular, for a convergent sequence with limiting eight-vertex probabilities $p_H$,
\begin{equation}
   \sum_H p_H T_{\sigma,K}(H)\ge0
   \qquad\text{for each term of }\Phi.
   \label{eq:certificate-limit}
\end{equation}

\subsection{The pointwise certificate}

Let $\mathcal P$ be the set of eight-vertex patterns arising in the two-part construction. Let $\mathcal T$ be the union of $\mathcal P$ and the patterns obtained from a two-part colouring of eight vertices by changing the colour of exactly one triple lying entirely in one part. These sets are defined up to relabelling.

\begin{lemma}[Finite certificate]\label{lem:certificate}
For the rational data in \texttt{certificate/proof.json}, the residual in~\eqref{eq:host-residual} vanishes on $\mathcal T$ and satisfies $R(H)\ge 10^{-5}$ for every $H\notin\mathcal T$. There are $893{,}924$ eight-vertex parity types, with $|\mathcal P|=22$ and $|\mathcal T|=72$.
\end{lemma}

\begin{proof}[Verification of the finite lemma]
An eight-vertex pattern is encoded by its 70 oddness bits. It satisfies the parity equation on each of the $\binom85=56$ five-sets. Fix a vertex $0$. The $\binom73=35$ bits $q(0ijk)$ determine the remaining bits, since the equation on $\{0,i,j,k,l\}$ expresses $q(ijkl)$ as the sum of the other four. Conversely any assignment of these 35 bits is realized by setting the triple bits through $0$ to zero and setting the colour of $ijk$ to $q(0ijk)$. The labelled universe thus has size $2^{35}$.

The certificate supplies one parity word for each proposed isomorphism class. To check this list, first test the 56 parity equations for each word. Treat its odd quartets as edges, and write $d_v$ and $d_{vw}$ for vertex degrees and pair codegrees. Give $v$ the signature consisting of $d_v$ and the sorted list of pairs $(d_w,d_{vw})$ over $w\ne v$. Sort the vertices by signature and permute within the equal-signature classes. For each resulting order record the 35-bit word through its first vertex, and take the least such word as a key. Isomorphic patterns give the same key. Conversely, a common key determines the same labelled parity pattern and hence implies isomorphism. The number of orders attaining the minimum is $|\operatorname{Aut}(H)|$, since the orders giving any fixed labelled pattern form an orbit under the free action of the automorphism group. The checker finds $893{,}924$ distinct keys and verifies
\[
       \sum_H\frac{8!}{|\operatorname{Aut}(H)|}=2^{35}.
\]
This checks exhaustion as well as the absence of duplicates.

For each word, the extension flags in~\eqref{eq:rooted-product} are read from the induced parity patterns. The arithmetic to be checked is a signed integer dot product with positive rational coefficients, which include the divisors in the table. A directed upper bound allows most comparisons to be made with word-sized integers. For an integer entry $M$, rational multiplier $a\ge0$, and scale $D$ one has
\[
       DaM\le \lfloor Da\rfloor M+\max(M,0).
\]
The verifier sums upper endpoints for $D\Phi(H)$ and compares the result with $\lfloor D(e(H)/70-5/16)\rfloor$. A strict inequality at this stage proves $R(H)>0$. The first implementation uses $D=2^{55}$ and reduces the unresolved set to 72 words. Multiplication by a common denominator gives exact equality in all 72 cases. The checker also derives $\mathcal P$ and $\mathcal T$ independently from the two-part rule by varying the part sizes and cross signs and, for each colouring, flipping each possible same-side triple. Canonicalization identifies the 72 zero residuals with precisely these 22 original and 50 additional patterns.

The second arithmetic implementation decodes the kernels afresh, checks triangle-freeness and uses $D=2^{49}$. For a six-set of $H$ form the graph whose edges are pairs complementary to odd quartets. Write $b_H$ for the fraction of six-sets on which this graph consists of three triangles sharing an edge, with an isolated vertex. The second implementation checks $R(H)\ge \eta\,\mathbf 1_{\{H\notin\mathcal T\}}(b_H+1/20)$, where $\eta>1/5000$ is the rational number \texttt{strict\_noncontact\_slack} in the certificate. This gives the asserted gap. The strengthened comparison requires exact rational arithmetic on 652 words. The implementations bound the magnitudes of integer accumulators before their use. The rational factor list and primitive parity words suffice for both computations.
\end{proof}

For a small example of the normalizations, suppose a colouring on eight vertices has just one red triple. Its five odd quartets all contain that triple. Among the 35 complementary quartet splits, five are odd--even and 30 are even--even. The unrooted term before weighting is therefore $[30(11)^2+5(11)(-5)]/35=671/7$. This pattern belongs to $\mathcal T$ by placing all eight vertices in one part and flipping the specified triple. Thus the full sum of weighted products is $65/70-5/16=69/112$ on this pattern. On the all-even pattern the corresponding value is $11/16$.

Averaging~\eqref{eq:host-residual} and applying~\eqref{eq:certificate-limit} yields, for every limiting colouring law,
\begin{equation}
       d-\frac5{16}\ \ge\sum_H p_H R(H)
           \ \ge\ 10^{-5}\sum_{H\notin\mathcal T}p_H .
       \label{eq:contact-gap}
\end{equation}
At equality the limiting law is supported on $\mathcal T$, and every individual rooted product with positive coefficient has limiting expectation zero.

\section{The finite-order estimate}\label{sec:finite}

The certificate evaluates a product on two disjoint extensions of a root.
Conditional positivity holds when the extensions are sampled independently.
To quantify the discrepancy for an ordinary square, it is useful to retain
the variance of its defining function.

\begin{lemma}\label{lem:disjoint-sampling}
Let $N\ge 2k\ge2$. Let $X$ be a uniform $k$-subset of $[N]$ and,
given $X$, let $Y$ be a uniform $k$-subset of $[N]\setminus X$.
For every $f\colon\binom{[N]}k\to\mathbb R$,
\[
 \mathbb E[f(X)f(Y)]
 \ \ge\ (\mathbb E f(X))^2-\frac{k}{N-k}\operatorname{Var}(f(X)).
\]
\end{lemma}

\begin{proof}
With the uniform inner product, the operator
$D\colon\mathbb R^{\binom{[N]}k}\to\mathbb R^{\binom{[N]}k}$ given by
$(Df)(X)=\mathbb E[f(Y)\mid X]$ is self-adjoint.
For $J\subseteq[N]$ put $f_J(X)=\mathbf 1_{J\subseteq X}$, and let $U_j$
be the span of the functions $f_J$ with $|J|\le j$. If $|J|=j\le k$, then
\[
 Df_J(X)=\frac{(k)_j}{(N-k)_j}\mathbf 1_{X\cap J=\varnothing}
        =\frac{(k)_j}{(N-k)_j}\sum_{I\subseteq J}(-1)^{|I|}f_I(X),
\]
where $(s)_j=s(s-1)\cdots(s-j+1)$. Thus $D$ preserves $U_j$ and acts on
$U_j/U_{j-1}$ by the scalar $(-1)^j(k)_j/(N-k)_j$.
Self-adjointness gives the same action on $U_j\cap U_{j-1}^{\perp}$.
The constants have eigenvalue one and $U_k$ is the full space. For $j\ge1$
the magnitude of the displayed ratio is at most $k/(N-k)$, because $N\ge2k$.
Apply this bound to the centered part of $f$.
\end{proof}

For a kernel $K$ coming from a triangle-free graph, write $a_i\ge0$ for the amplitude at
extension state $i$. Its diagonal entry is $a_i^2$, and
$|K(i,j)|\le (a_i^2+a_j^2)/2$. Given a root, first take two independent
uniform outside $k$-sets. Their expected kernel value is nonnegative.
Conditional on an overlap, either outside set is still marginally uniform.
Consequently the expected kernel on a disjoint pair is bounded below by
minus the mean diagonal entry times the ratio of the probabilities of overlap
and disjointness. For $k=1$ this ratio is $1/(N-1)$. For $k=2$ it is
\[
      \frac{\binom N2-\binom{N-2}2}{\binom{N-2}2}
        =\frac{4N-6}{(N-2)(N-3)}\le\frac5{N-3}\qquad(N\ge4).
\]

\begin{proposition}\label{prop:finite-collision}
For every colouring $G$ on $n\ge8$ vertices, $d(G)\ge5/16-L_n$, where
\[
 L_n=
 \frac{0.216198}{n-4}+\frac{1.662165}{n-5}
 +\frac{2.708593}{n-6}+\frac{1.175666}{n-7}.
\]
\end{proposition}

\begin{proof}
For terms with root size $r\in\{0,2,4,6\}$ the number of outside
vertices is $N=n-r$ and the extension size is $k=(8-r)/2$.
For each group of terms in Table~\ref{tab:finite-diagonals}, form the sum
of its weighted diagonal entries on a single $(r+k)$-vertex flag. Average
over the root and its orders inside a fixed $(r+k)$-vertex pattern, with the
normalization of the certificate. The maximum over such patterns is denoted by
$M_r$ for the ordinary terms at sizes $0,2,4$, by $M_{4,c}$ for the
triangle-free terms at size four, and by $M_6$ for all terms at size six.
For $M_0$ one can simply use $121$ times the coefficient of the unrooted
contrast. All these maxima bound the corresponding mean diagonals for every
finite colouring.

\begin{table}[htbp]
\centering
\begin{tabular}{lll}
\hline
Source of error & Denominator & Rational upper bound on numerator\\
\hline
Unrooted squares                     & $n-4$ & $4M_0\le216198/10^6$\\
Two-root squares                     & $n-5$ & $3M_2\le1662165/10^6$\\
Four-root ordinary squares           & $n-6$ & $2M_4\le2708593/10^6$\\
Six-root terms and four-root kernels & $n-7$ & $M_6+5M_{4,c}\le1175666/10^6$\\
\hline
\end{tabular}
\caption{Diagonal sums in the finite-order estimate. Each upper bound
is rounded upwards from an exact rational calculation.}
\label{tab:finite-diagonals}
\end{table}

For the two-root terms the calculation averages over the $20$ ordered pairs
of a five-set. Its quartet bits have weight $0$, $2$ or $4$.
For the four-root terms it suffices to check the $16$ six-vertex types.
For the six-root terms the verifier enumerates the $2^{20}$ possible labelled
seven-vertex parity words and averages over the seven choices of the outside
vertex. The exact rational sums are supplied with the certificate.
Lemma~\ref{lem:disjoint-sampling} gives the first three errors in the table
and the error for ordinary six-root squares. The overlap ratios above give
the other errors. Thus $\mathbb E\Phi(H)\ge-L_n$ on a uniformly sampled
eight-set $H$. Averaging the pointwise certificate proves the proposition.
\end{proof}

\section{The law at equality}\label{sec:equality}

We first determine the support of an extremal law on eight vertices. An all-even five-set will then allow the local bipartitions on larger samples to be recovered from their quartet parities.

\subsection{Removing the fifty additional types}
At equality the certificate puts the limiting eight-vertex law on $\mathcal T$. To restrict it to $\mathcal P$, we use one more indicator on a six-set with four of its vertices distinguished. Form the graph on the six-set by joining two vertices when the complementary quartet is odd. The indicator $F$ is one precisely when its degree multiset and the degrees of the four root vertices occur on the same row of the following table. The root degrees are measured in the full six-vertex graph.
\[
\begin{array}{c|l|c}
\text{six vertex degrees}&\text{four root degrees}&
    \text{admissible roots out of 15}\\\hline
0,0,0,2,2,2 & 0,0,2,2 &9\\
0,2,2,2,4,4 & 0,2,2,4 &6\\
2,2,2,4,4,4 & 2,2,2,4\ \text{or}\ 2,2,4,4 &12
\end{array}
\]
The three rows describe, respectively, a triangle with three isolates, the three-page book with one isolate, and a six-cycle with a triangle on three alternate vertices. For the last identification, the subgraph on the three degree-four vertices must have three more edges than that on the other three vertices. Thus it is a triangle, the other three form an independent set, and the edges between them make a six-cycle.

The root choices reflect how a single pure flip can create one of these graphs. Before a flip the possible edge counts, indexed by the smaller part size $0,1,2,3$, are $0$, $10$, $8,11,12$, and $6,10,12$, respectively. Thus all three displayed graphs require a flip within the six-set. If three of the six vertices are on the side of the flipped triple, they are exactly that triple. After the flip they span a triangle, the other three vertices are independent, and the cross edges have even row and column sums. The cross graph is empty, a four-cycle, or a six-cycle. In the three cases the flipped triple consists of the triangle, both book centers and one leaf, or all three degree-four vertices, respectively. With four vertices on the side of the flip there are initially eight cross edges. The flip removes two and adds an edge on the opposite side. Any additional edges among the four vertices number zero, three, or four, according to the cut given by the products of the two cross signs in each row (the complementary pair indexes the edge). Only the book can occur among the three displayed graphs, and the flip is its three leaves. With five vertices on the side of the flip the result has eleven edges. With six vertices on that side, the result is a triangle and the flip consists of its three isolates. Inspection of the table now shows that none of the admissible roots contains a possible flipped triple. This makes it impossible for one flip to produce two marked copies whose intersection is their common root.

For an eight-type $H$, let $s(H)$ be the proportion of the $70\cdot6=420$ ways to choose a four-root and a first outside pair in which $F=1$. Let $t(H)$ be the proportion of the $70\cdot3=210$ ways to choose a root and split its complement into pairs in which both pairs give $F=1$.

\begin{lemma}\label{lem:eight-law}
In any limiting parity law, put $\varepsilon=\sum_{H\notin\mathcal T}p_H$. Then
\[
                  \sum_{H\notin\mathcal P}p_H
                         \le \varepsilon+70\sqrt{\varepsilon}.
\]
In particular every eight-vertex sample from an equality law belongs to $\mathcal P$ almost surely.
\end{lemma}

\begin{proof}
For this indicator, direct enumeration gives
\[
\begin{array}{c|cc}
H &s(H)&t(H)\\\hline
H\in\mathcal P & 0 &0\\
H\in\mathcal T\setminus\mathcal P &\ge 6/420 &0.
\end{array}
\]
The calculation uses only the $72$ words of $\mathcal T$. On each of their $28$ six-subsets read the $15$ quartet bits as complementary edges, compute the six degrees, and mark the roots in the table. For each four-root count the marked pairs, then check the three splits of its complement. The minimum positive number of marked root--pair choices is six. None of the $72\cdot210=15{,}120$ splits has both pairs marked. The $22$ words with no marks are exactly the two-part words, as one checks directly by trying each bipartition and fixing the signs in one row and column of its cross matrix.

In a finite colouring, let $m(R)$ be the fraction of outside pairs marked at a uniform four-root $R$. Then $\mathbb E_R m(R)=\mathbb E s(H)$ for a uniform eight-set $H$. For two independent outside pairs, the probability of two marks is $\mathbb E_R m(R)^2$. Requiring these pairs to be disjoint changes that probability by at most
\[
                  1-\binom{n-6}{2}\big/\binom{n-4}{2}.
\]
Thus in the limit $(\mathbb E s(H))^2\le\mathbb E t(H)$ by Jensen's inequality. The finite calculation gives $\mathbb E t(H)\le\varepsilon$ and $\Pr[H\in\mathcal T\setminus\mathcal P]\le70\mathbb E s(H)$. Together these prove the bound. At equality, $\varepsilon=0$ by~\eqref{eq:contact-gap}.
\end{proof}

A finite \emph{two-part representation} consists of a partition $A\cup B$ and signs $Z_{ab}$, $a\in A$, $b\in B$. Pure quartets are even, quartets split three to one are odd, and on $i,j\in A$, $x,y\in B$ one has $(-1)^{q(ijxy)}=Z_{ix}Z_{iy}Z_{jx}Z_{jy}$. Empty parts are allowed. Every such representation on eight vertices occurs with positive probability in the fair construction. Hence $\mathcal P$ is precisely its eight-vertex support.

\subsection{Recovering the parts and their cross signs}
\begin{lemma}\label{lem:glue}
If every eight-set of a sufficiently large finite parity pattern has type in $\mathcal P$, then the whole pattern has a two-part representation.
\end{lemma}

\begin{proof}
The finite Ramsey theorem \cite{Ramsey}, applied to the two colours of the quartets, supplies a five-set $P$ on which all quartet bits agree. The five-set parity identity forces the common bit to be zero. We will use only the existence of this one all-even five-set.

In every representation containing $P$, its vertices lie on a common side, since any five-set meeting both parts contains a quartet split three to one. Call the side of $P$ the $A$ side. For $v\notin P$, the parities $q(ijkv)$ for triples $i,j,k\in P$ are all equal, by extending $P\cup\{v\}$ to an eight-set and using its representation. Declare $v\in A$ when this common bit is zero, and $v\in B$ otherwise. In any representation containing $P$, this rule recovers the sides of every other vertex in that representation.

We first check the pure and three-to-one rules. A quartet meeting $P$ can be checked in a representation containing $P$ and the quartet. For a quartet $S$ disjoint from $P$ that contains some $v\in A$, take four vertices of $P$ together with $S$. In a representation of this eight-set, the four chosen vertices together with $v$ form an all-even five-set, so they lie on a common side. The three-anchor test just established puts every other vertex of $S$ on its prescribed side. This settles all such quartets. If $S\subseteq B$, adjoin one vertex of $P$. The four quartets containing it are odd by the preceding tests, so the fifth quartet $S$ is even.

It remains to recover a single cross matrix for all two-to-two quartets. Work in $\mathbb F_2$. For distinct $i,j,k\in A$ and $x,y\in B$, the five-set identity gives
\[
                 q(ijxy)+q(jkxy)+q(kixy)=0,
\]
since the other two quartets are odd. The same identity holds with the two parts exchanged. If both parts are nonempty, fix $a_0\in A$ and $b_0\in B$. Set $z_{a_0x}=z_{ib_0}=0$ and $z_{ix}=q(a_0ib_0x)$ for the remaining entries. Applying the two identities, first in the $A$ pair and then in the $B$ pair, gives
\[
                    q(ijxy)=z_{ix}+z_{iy}+z_{jx}+z_{jy}.
\]
The equations with an anchor as endpoint follow directly by the same identity in the other pair. The signs $Z_{ix}=(-1)^{z_{ix}}$ give the representation. If a part is empty the pure rule suffices.
\end{proof}

\begin{proof}[Proof of Theorem~\ref{thm:main-law}]
Extract a subsequence with convergent labelled sample distributions at every finite order, and write $Q_k$ for a random pattern from the limiting $k$-vertex distribution. These distributions are consistent under uniform restriction. By Lemma~\ref{lem:eight-law} and a union bound, every eight-subset of $Q_k$ has type in $\mathcal P$ almost surely. For all sufficiently large $k$, choose a two-part representation by Lemma~\ref{lem:glue}.

Let $\alpha,\beta$ be the part fractions, let $h=\alpha-\beta$, and let
\[
 t=\mathbb E_{a,a'\in A}\Bigl(\mathbb E_{b\in B}Z_{ab}Z_{a'b}\Bigr)^2
\]
with independent uniform choices and $t=0$ if a part is empty. This nonnegative quantity is the mean signed rectangle when vertices are sampled with replacement. Sampling four vertices with replacement incurs a collision with probability at most $6/k$. Counting by part sizes therefore gives, uniformly in the chosen representation,
\begin{align*}
 d(Q_k)
  &=\alpha^4+\beta^4+3\alpha^2\beta^2(1+t)+O(k^{-1})\\
  &=\frac5{16}+\frac{6h^2+5h^4}{16}+3\alpha^2\beta^2t+O(k^{-1}).
\end{align*}
Consistency gives $\mathbb E d(Q_k)=5/16$. Hence $h\to0$ in probability. On the event $|h|<1/2$, the coefficient of $t$ is bounded away from zero, so the same identity gives $t\to0$ in probability.

To see what this implies for the cross signs, let $f\colon A\to[-1,1]$ and $g\colon B\to[-1,1]$. Two applications of Cauchy--Schwarz give
\[
                    |\mathbb E_{a,b}f(a)Z_{ab}g(b)|\le t^{1/4}.
\]
For any fixed nonempty simple bipartite edge set, sample its vertices independently within their prescribed parts and take the product of the corresponding cross signs. Condition on every vertex except the endpoints of one selected edge. All other factors involving those endpoints separate into a function of each, so the same $t^{1/4}$ bounds the absolute mean of the product. Thus every nonconstant sign moment vanishes in the limit. Expanding the indicator of a particular sign assignment in these moments shows that the sampled cross signs tend jointly to independent fair signs. Sampling instead without replacement costs at most $\binom r2/k$ for an $r$-vertex sample from all $k$ vertices. Its part labels likewise tend to independent fair labels because $h\to0$.

By consistency this identifies the law of $Q_r$ itself. The argument applies to every subsequential limit and proves the claimed convergence.
\end{proof}

\paragraph{An independent determination of the eight-vertex probabilities.}
The equality calculation can also be made directly on $\mathcal T$. For an ordinary four-root square with extension function $u$, let $m_u(R)$ be its conditional mean on an outside pair, taking $u=0$ at roots of the wrong type. At equality the square has vanishing expectation, so $\mathbb E_R m_u(R)^2\to0$. Multiplication by the conditional probability of any compatible extension flag $j$ consequently gives a vanishing mixed moment. Its coefficient at $H$, after scaling by $8!$, is the sum of $u$ on the first extension times the indicator of $j$ on the second over all vertex orders of $H$. The first four positions specify the ordered root and the final four positions the two pairs.

The supplementary computation selects $70$ such rows from five positively weighted squares. Two further rows record total mass and density, with entries $1$ and $16e(H)-350$. In the supplied ordering their $72$-by-$72$ integer matrix has determinant $757006\bmod 1000003$. Thus the probabilities are unique already at order eight. For comparison the fair model gives the probability of a labelled word by summing over its compatible ordered splits of sizes $a,b$. A split with $a,b>0$ contributes $2^{-8-(a-1)(b-1)}$, since fixing the first row and column leaves $(a-1)(b-1)$ independent rectangle bits. An empty part contributes $2^{-8}$ if the word is all-even. Multiplication by the size of the labelled orbit gives the unlabelled probability. These probabilities have common denominator $65536$ and satisfy the integer system exactly.

\section{The first finite-size correction}\label{sec:second-order}

The four bounds in Table~\ref{tab:finite-diagonals} apply to every
colouring. At a minimizing colouring we also know the limiting law of
the flags responsible for the sampling error. This gives a much narrower
interval for the coefficient of $1/n$. Put $c=5/16$.

\begin{proposition}\label{prop:second-order}
The minimum even-quartet density satisfies
\[
 \frac34\ \le\ \liminf_{n\to\infty} n(c-m_n)
 \ \le\ \limsup_{n\to\infty} n(c-m_n)\ \le\ 0.752723 .
\]
\end{proposition}

\begin{proof}
For the first inequality, take any Hadamard matrix of order $M$ as the
cross-sign matrix in the construction with two parts of size $M$.
Counting rectangles on pairs of orthogonal rows gives
\[
   c-d(G_{2M})=\frac{3(8M-11)}{16(2M-1)(2M-3)}
               =\frac{3/4+o(1)}{2M}.
\]
There are Hadamard orders $M\ge n/2$ with $2M/n\to1$. Here is an
elementary verification. Order $2$ is immediate. For order $12$ use the Paley construction~\cite{Paley}. Put $Q_{xy}=\chi(y-x)$ for $x,y\in\mathbb F_{11}$, with $\chi\colon\mathbb F_{11}\to\{-1,0,1\}$ the
quadratic character and $\chi(0)=0$. Then $Q+Q^{\mathsf T}=0$,
$Q\mathbf 1=0$ and $QQ^{\mathsf T}=11I-J$. For distinct $x,y$, the last
identity uses the substitution $t\mapsto(t-x)/(t-y)$ from
$\mathbb F_{11}\setminus\{y\}$ onto $\mathbb F_{11}\setminus\{1\}$.
Adjoining a first row and column of ones to $Q-I$ gives the desired
matrix. Tensor products give all orders $2^a12^b$ for nonnegative $a,b$.
Since $\log_2 12$ is irrational, for every $\varepsilon>0$ finitely many
multiples $b\log_2 12$ meet every interval of length $\varepsilon$ modulo
one. Thus for all sufficiently large $x$ there are $a,b\ge0$ with
$x\le a+b\log_2 12\le x+\varepsilon$. Apply this with $x=\log_2(n/2)$
and use $m_n\le d(G_{2M})$ by uniform deletion.

For the upper bound take a minimizing colouring at each $n$. Consider
one term of $\Phi$ with root size $r$ and extension size $k=(8-r)/2$. Include its weight, root indicator and factor $\kappa_r$ in the value being averaged. Sample the root as in the
certificate and choose its two outside $k$-sets independently. Let $A_j$
be the mean kernel value conditional on an intersection of size $j$, and
$p_j$ the probability of that intersection. The certificate samples the
term with expectation $A_0$, whereas $\sum_jp_j A_j\ge0$ by conditional
positivity. All $A_j$ are bounded independently of $n$. For $N=n-r$,
\[
   \frac{p_1}{p_0}
     =\frac{k\binom{N-k}{k-1}}{\binom{N-k}{k}}
     =\frac{k^2}{n-7},
   \qquad \sum_{j\ge2}\frac{p_j}{p_0}=O(n^{-2}).
\]
Hence $A_0\ge-k^2 A_1/(n-7)-O(n^{-2})$. This estimate retains the
value on a single overlap, including its sign.

A single overlap uses exactly seven vertices. By Theorem~\ref{thm:main-law},
its limiting expectation at these minimizers is the corresponding
expectation $A_1^*$ in the balanced model. Sum the last inequality over
the $841$ terms and use $R(H)\ge0$. It follows that
\[
 \limsup_{n\to\infty} n(c-m_n)\le C_*,
 \qquad C_*:=\sum_{\mathrm{terms}} k^2 A_1^* .
\]
For precision, the entries in the following display are rational
\emph{upper bounds} on the four subtotals, obtained by rounding each up
at the indicated precision.
\[
\begin{array}{c|rrrr}
r&0&2&4&6\\ \hline
\text{upper bound}&0&0.093837982414&0.442961486423&0.215923418603
\end{array}
\]
To evaluate $C_*$, take each of the $22$ eight-vertex patterns with its
known probability, average over its $8!$ orders, and ignore the last
vertex. Use the first $r$ vertices as root, the next as the common
outside vertex and two consecutive disjoint groups of $k-1$ remaining
vertices to complete the extensions. Apply the same weights and
normalizations as in the certificate. Thus these values require only
integer counts of the original flag pairs on seven labelled vertices.
The exact arithmetic is reproduced by \texttt{collision\_limit\_checker.py}.
An independent calculation enumerates all $128$ labelled bipartitions on
seven vertices, sums over the allowed cross matrices after fixing one row
and column, and obtains the same four rational values.
The displayed upper bounds sum to less than $0.752723$.
\end{proof}

\section{A quadratic bound in the six-vertex book density}\label{sec:book}

On a six-set $S$ form the graph $D(S)$ in which $xy$ is an edge exactly when
$q(S\setminus\{x,y\})=1$. Let $B$ be the graph consisting of three triangles
with a common edge together with an isolated vertex. Its two vertices of
degree four will be called the centers and its three vertices of degree two
the leaves. Write $b$ for the limiting density of six-sets $S$ with
$D(S)\cong B$.

\begin{proposition}\label{prop:quadratic-book}
For every limiting pair $(d,b)$ of even-quartet and book densities,
\[
             d\ \ge\ \frac5{16}+\frac{144117}{1074218750}\,b^2 .
\]
\end{proposition}

\begin{proof}
Distinguish a four-set $R\subset S$. Call $(R,S)$ a rooted book if $D(S)$
is a book and $R$ consists of its isolated vertex, one center and two leaves.
There are $2\binom32=6$ such roots out of $15$. Every such root is odd, since
the complementary pair in $D(S)$ consists of a center and a leaf.

On an eight-type $H$, let $e_F(H)$ be the fraction of the $70\cdot3$ choices
of a four-set $R$ and an unordered split of its complement into pairs $P,Q$
for which both $(R,R\cup P)$ and $(R,R\cup Q)$ are rooted books. The
supplementary integer check proves
\[
        R(H)\ \ge\ \lambda e_F(H),\qquad
        \lambda=\frac{144117}{250000000}.
\]
Here $R(H)$ is the residual in Lemma~\ref{lem:certificate}.
The check constructs the indicator from the six-vertex quartet bits and
counts successful splits in each eight-type. Exactly $135{,}740$ types
have positive count. On these types the directed integer lower bound on
the residual suffices for the displayed comparison. On all remaining types
it follows from $R(H)\ge0$.

In a colouring on $n$ vertices, condition on a uniform four-root $R$ and
let $p_R$ be the probability that a uniform outside pair completes a rooted
book. If the finite densities are $d_n$ and $b_n$, then
$\mathbb E p_R=2b_n/5$ and $p_R$ vanishes on even roots. Thus
\[
        (2b_n/5)^2\le (1-d_n)\,\mathbb E p_R^2.
\]
Two independent outside pairs overlap with probability $O(1/n)$, uniformly
in $R$. It follows that $\mathbb E e_F(H)=\mathbb E p_R^2+O(1/n)$.
Since $\mathbb E\Phi(H)\ge-L_n$, passage to a limiting law with $d<1$
gives
\[
       d-\tfrac5{16}\ge\lambda\,\frac{4b^2}{25(1-d)}
                  \ge\frac{64\lambda}{275}\,b^2 .
\]
The last step uses $d\ge5/16$. When $d=1$, the preceding
Cauchy--Schwarz inequality gives $b=0$ on taking the limit.
\end{proof}

\begin{corollary}\label{cor:book-power}
The exponent two in Proposition~\ref{prop:quadratic-book} is optimal. For every $0<p<2$ and every $C>0$ there is a limiting law with $b>0$ and $d-5/16<Cb^p$.
\end{corollary}

\begin{figure}[htb]
\centering
\begin{minipage}[c]{.28\linewidth}\centering
$\begin{array}{c|rrr}
 & b_1&b_2&b_3\\ \hline
 a_1&-&+&+\\
 a_2&+&+&+\\
 a_3&+&+&+
\end{array}$
\end{minipage}
\hspace{.06\linewidth}
\begin{minipage}[c]{.47\linewidth}\centering
\begin{tikzpicture}[x=1cm,y=1cm,vertex/.style={circle,fill,inner sep=1.8pt},every label/.style={font=\small}]
\coordinate (u) at (0,0); \coordinate (v) at (1.75,0);
\coordinate (l) at (.875,.66); \coordinate (m) at (.875,1.37); \coordinate (n) at (.875,-.70);
\coordinate (i) at (2.7,.6);
\draw (u)--(v) (u)--(l)--(v) (u)--(m)--(v) (u)--(n)--(v);
\node[vertex,label=left:$a_2$] at (u) {};
\node[vertex,label=right:$a_3$] at (v) {};
\node[vertex,label=below:$a_1$] at (l) {};
\node[vertex,label=above:$b_2$] at (m) {};
\node[vertex,label=below:$b_3$] at (n) {};
\node[vertex,label=above:$b_1$] at (i) {};
\end{tikzpicture}
\end{minipage}
\caption{A single flip of $a_1a_2a_3$ produces the complementary graph on the right. Its common edge is $a_2a_3$ and $b_1$ is isolated.}
\label{fig:book}
\end{figure}
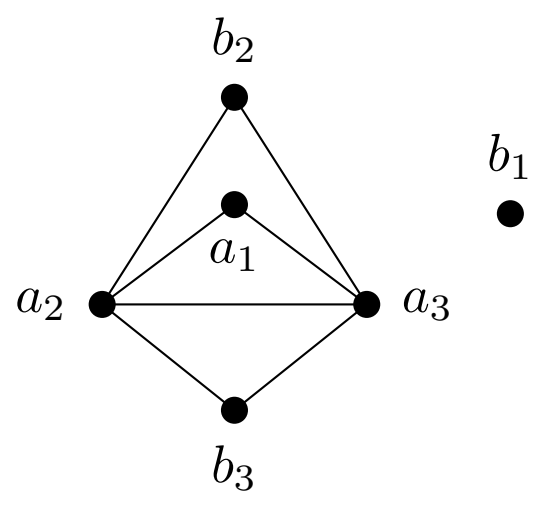

\begin{proof}
Use the perturbation of Section~\ref{sec:perturbation}. On a fixed six-set consider the event that three vertices lie in each part, and order the vertices in each part by their original labels. Prescribe a cross matrix whose only negative entry is $Z_{a_1b_1}$. This has conditional probability $2^{-9}$. The unique $A$-triple flips with probability $t$. In the complementary graph the three $A$-pairs are edges, the three $B$-pairs lose their edges, and $a_i b_j$ is an edge exactly when $i,j\in\{2,3\}$. This is the book in Figure~\ref{fig:book}. Thus $b(t)\ge\binom63 2^{-6}2^{-9}t$. Since $d(t)-5/16=O(t^2)$ by~\eqref{eq:perturb}, the ratio $(d(t)-5/16)/b(t)^p$ tends to zero.
\end{proof}

\section{Verification and further questions}\label{sec:verification}

The electronic supplement contains the rational coefficients, the parity words representing the eight-vertex types, source code, and an entry point for repeating the finite calculations. From the extracted archive run
\begin{center}
\texttt{python3 verify\_all.py}
\end{center}
on a little-endian x86-64 Linux system with Python 3.11 or later, NumPy, and a C++17 compiler available as \texttt{c++}. No optimization solver is needed. Lemma~\ref{lem:certificate} describes the orbit test and the arithmetic comparisons. These checks certify the supplied coefficients irrespective of how they were found. The remaining programs check the rooted six-vertex indicator, the book inequality, the finite-order constants, and the first finite correction. The supplementary nine-vertex census checks every one-vertex extension of each of the $22$ types and finds an all-even five-set in every admissible extension. Thus the anchor proof of Lemma~\ref{lem:glue} applies already on nine or more vertices. The Ramsey argument there suffices for Theorem~\ref{thm:main-law}. The rational witness \texttt{certificate/proof.json} has SHA-256 digest
\begin{center}\small
\texttt{5e6b2f8e662abbc6b22ee9202efe5860c2f98879d831586a77761b95c63e237a}.
\end{center}

The bounds in Proposition~\ref{prop:second-order} leave a short interval for the finite correction. Is $m_n=5/16-3/(4n)+o(1/n)$? Orthogonal sign matrices attain the proposed correction asymptotically, and equality in the limiting problem determines every fixed sample distribution. The question is whether the loss from sampling disjoint extensions can be bounded at first order by the amount seen in those matrices. For unrestricted coverings, even the leading coefficient would have to be established without the five-set parity identity.

\section*{Acknowledgments}
The author thanks OpenAI for the AI tools used in developing and verifying the computational certificates and in preparing this manuscript.

\end{document}